\documentclass[11pt]{article}

\usepackage{graphicx,times,latexsym,amssymb,amsbsy,amsmath,epsfig}

\input xy

\xyoption{all}

\newtheorem{thm}{Theorem}[section]
\newtheorem{prop}[thm]{Proposition}
\newtheorem{lem}[thm]{Lemma}

\newtheorem{dfn}[thm]{Definition}

\newtheorem{ques}[thm]{Question}
\newtheorem{conj}[thm]{Conjecture}

\newcommand{\bs}[1]{\boldsymbol{#1}}
\renewcommand{\bf}[1]{\mathbf{#1}}
\renewcommand{\rm}[1]{\mathrm{#1}}

\newcommand{\bbN}{\mathbb{N}}

\newcommand{\bbQ}{\mathbb{Q}}
\newcommand{\bbR}{\mathbb{R}}

\newcommand{\bbZ}{\mathbb{Z}}

\newcommand{\bphi}{\boldsymbol{\phi}}
\newcommand{\bpsi}{\boldsymbol{\psi}}

\newcommand{\sfC}{\mathsf{C}}
\newcommand{\sfE}{\mathsf{E}}

\newcommand{\sfZ}{\mathsf{Z}}

\newcommand{\rmF}{\mathrm{F}}

\renewcommand{\d}{\mathrm{d}}

\newcommand{\id}{\mathrm{id}}

\newcommand{\frH}{\mathfrak{H}}

\newcommand{\G}{\Gamma}

\renewcommand{\a}{\alpha}

\newcommand{\eps}{\varepsilon}

\newcommand{\s}{\sigma}

\newcommand{\uhr}{\upharpoonright}

\newcommand{\qed}{\nolinebreak\hspace{\stretch{1}}$\Box$}

\newcommand{\actson}{\curvearrowright}
\newcommand{\barint}{-\!\!\!\!\!\!\int}
\renewcommand{\vec}[1]{\overrightarrow{#1}}
\renewcommand{\t}[1]{\tilde{#1}}

\begin{document}

\title{Norm convergence of continuous-time polynomial multiple ergodic averages}

\author{TIM AUSTIN\\ \small{Brown University}\\ \small{Providence, RI, USA}\\ \small{\texttt{timaustin@math.brown.edu}}\\ \\ \normalsize{\emph{Dedicated to the memory of Dan Rudolph}}}

\date{}

\maketitle

\parindent 0pt

\begin{abstract} For a jointly measurable
probability-preserving action $\tau:\bbR^D\curvearrowright (X,\mu)$
and a tuple of polynomial maps $\bf{p}_i:\bbR\to \bbR^D$,
$i=1,2,\ldots,k$, the multiple ergodic averages
\[\frac{1}{T}\int_0^T (f_1\circ \tau^{\bf{p}_1(t)})(f_2\circ\tau^{\bf{p}_2(t)})\cdots (f_k\circ\tau^{\bf{p}_k(t)})\,\d t\]
converge in $L^2(\mu)$ as $T\to\infty$ for any $f_1,f_2,\ldots,f_k
\in L^\infty(\mu)$.  This confirms the continuous-time analog of the conjectured norm convergence of discrete polynomial multiple ergodic averages, which in is its original
formulation remains open in most cases.  A proof of convergence
can be given based on the idea of passing up to a sated extension of $(X,\mu,\tau)$ in
order to find a simple partially characteristic factor, similarly to the recent development of this idea for the study of related discrete-time averages, together with a new
inductive scheme on tuples of polynomials.  The new induction scheme becomes available upon changing the time variable in the above integral by some fractional power, and provides an
alternative to Bergelson's PET induction, which has been the
mainstay of positive results in this area in the past.
\end{abstract}


\parskip 7pt

\section{Introduction}

Given commuting probability-preserving transformations $T_1,T_2,\ldots,T_k
\curvearrowright (X,\mu)$, the study of the associated `multiple'
(or `diagonal') ergodic averages
\[\frac{1}{N}\sum_{n=1}^N (f_1\circ T_1^n)(f_2\circ T_2^n)\cdots (f_k\circ T_k^n)\]
by now has an extensive history.  Interest in them originated in
Furstenberg's proof in~\cite{Fur77} of the Multiple Recurrence
Theorem, which effectively shows that the integrals of the above averages
stay uniformly positive as $N\to\infty$ when $T_i  = T^i$ for some
fixed $T$ and $f_1 = f_2 = \ldots = f_k = 1_A$ for some
non-negligible set $A$.  That work was followed by a
multidimensional generalization in Furstenberg and Katznelson's
paper~\cite{FurKat78}, and since then the above averages, various
generalizations of them and a host of related questions have been
investigated by several researchers: see, in particular, the papers
\cite{ConLes84,ConLes88.1,ConLes88.2}, \cite{Ber87}, \cite{Rud93},
\cite{Zha96}, \cite{FurWei96}, \cite{BerLei02},
\cite{HosKra01,HosKra05,HosKra05poly}, \cite{Zie07}, \cite{FraKra05}
and \cite{Tao08(nonconv)}.

In a sequence of breakthroughs culminating in Tao's work, it has now
been shown that the above averages always converge in $L^2(\mu)$ as
$N\to\infty$, augmenting our understanding of their asymptotic
nonvanishing due to Furstenberg and Katznelson.  However, their
generalization to `polynomial' multiple averages, such as
\[\frac{1}{N}\sum_{n=1}^N (f_1\circ T_1^{p_1(n)})(f_2\circ T_2^{p_2(n)})\cdots (f_k\circ T_k^{p_k(n)})\]
for $p_1,p_2,\ldots,p_k \in \bbZ[X]$, remain far less well
understood.  In case $T_1 = T_2 = T_3 = \ldots = T_k$ their
convergence and also a reasonably good description of their limit
have been obtained in~\cite{Lei05(poly),HosKra05poly}, and similarly in the case of more general commuting transformations under some assumptions of weak mixing~\cite{Ber87} or of some algebraic independence among the polynomials~\cite{ChuFraHos09}.  A handful of
other higher-dimensional cases are now at least partly
understood~\cite{Chu09,Aus--lindeppleasant1,Aus--lindeppleasant2}.
However, the following broader conjecture remains open in general:

\begin{conj}\label{conj}
For any probability-preserving action $T:\bbZ^D \curvearrowright
(X,\mu)$ and any polynomial mappings
$\bf{p}_1,\bf{p}_2,\ldots,\bf{p}_k:\bbZ\to \bbZ^D$ the averages
\[\frac{1}{N}\sum_{n=1}^N (f_1\circ T^{\bf{p}_1(n)})(f_2\circ T^{\bf{p}_2(n)})\cdots (f_k\circ T^{\bf{p}_k(n)})\]
converge in $L^2(\mu)$ for any $f_1,f_2,\ldots,f_k \in
L^\infty(\mu)$.
\end{conj}

A version of this appears as Question 9 in Bergelson~\cite{Ber96}, although it had certainly been proposed informally before that, by Furstenberg and others.  In addition to the special cases mentioned above, it has been proved for general polynomial mappings under the assumption that the action $T$ is totally ergodic by Johnson in~\cite{Joh09}.  In~\cite{BerLei02}, Bergelson and Leibman extend this conjecture even further
to the setting of actions $\G\curvearrowright (X,\mu)$ of a discrete
nilpotent group $\G$ and `polynomial mappings'
$\bf{p}_i:\bbZ\to\G$; here we will omit the technical preliminaries
needed to set up this latter notion.

This conjecture has been verified in all special cases that have
been successfully analyzed to date; but the partial results obtained
in~\cite{Aus--lindeppleasant2} indicate that the general conjecture may still lie some way beyond
the scope of current approaches to these results.  However, in this
paper we will see that the `continuous-time' analogs of the above
averages are rather simpler to understand, and do indeed enjoy the
analogous convergence.

\begin{thm}\label{thm:main}
For any jointly measurable action $\tau:\bbR^D \curvearrowright (X,\mu)$ and any tuple
of polynomials $\bf{p}_1$, $\bf{p}_2$, \ldots,
$\bf{p}_k:\bbR\to\bbR^D$ the associated multiple ergodic averages
\[A_T(f_1,f_2,\ldots,f_k) := \barint_0^T \prod_{i=1}^k f_i\circ \tau^{\bf{p}_i(t)}\,\d t\]
converge in $L^2(\mu)$ as $T\to\infty$ for any $f_1$, $f_2$, \ldots,
$f_k \in L^\infty(\mu)$, where $\barint_I$ denotes the average $\frac{1}{|I|}\int_I$ for any bounded real interval $I$.
\end{thm}

By analogy with Bergelson and Leibman~\cite{BerLei02}, it should be
possible to formulate an extension of the above conclusion that
applies to a tuple of polynomial maps $\bf{p}_i:\bbR\to N$ taking
values in some nilpotent Lie group $N$ and an action
$N\curvearrowright (X,\mu)$; however, this seems to lead to new
algebraic complications that we will not try to surmount in this
paper.

Theorem~\ref{thm:main} will be proved by induction on the set of polynomial maps
$\{\bf{p}_1,\bf{p}_2,\ldots,\bf{p}_k\}$.  An important invention
from~\cite{Ber87} is Bergelson's `Polynomial Ergodic Theorem' (`PET')
induction scheme: a wellordering on such sets of polynomials that he
and many others have now used as the basis for analyzing multiple
polynomial ergodic averages.  However, in this paper we depart from
this scheme, using instead a simple but crucial change in the
time-variable that allows us to work instead with families of
`fractional polynomial maps', among which a different and rather
more efficient induction scheme becomes more natural.  In a sense, the availability of these fractional polynomial maps in the continuous-time setting indicates a much greater `smoothness' that is enjoyed by the averages of Theorem~\ref{thm:main} than by their discrete-time analogs, hence leading to simplified behaviour.  This flexibility in choosing time-changes was also shown to me by Vitaly Bergelson.

In Section~\ref{sec:right-ques} we will introduce and justify these changes of
variables, and reformulate Theorem~\ref{thm:main} accordingly.  Then, in Section 3, we introduce an extension of that reformulated result which will be needed in order that our induction close on itself, and then we introduce the partial ordering on families of (fractional) polynomials that we will use and complete the proof.  In addition to the time-change and this partial order, the main tool we need is the technology of sated extensions developed in~\cite{Aus--lindeppleasant1}, which in our setting enables us to extend an initially-given $\bbR^D$-system to one in which the averages of interest exhibit
simplified behaviour from which their convergence can be deduced.

Before launching into technical details, we should mention two other recent works concerning multiple averages in continuous-time.  The first, that of Potts~\cite{Pot09}, studies the averages
\[\barint_0^T \prod_{i=1}^k (f_i\circ \tau^{it})\,\d t\]
for a single flow $\tau:\bbR\actson (X,\mu)$.  She shows that for these averages the powerful machinery developed by Host and Kra (see~\cite{HosKra05}, and also~\cite{HosKra05poly}) on multiple averages for a single transformation has a direct counterpart (and in fact that the Host-Kra-like factors of the system $(X,\mu,\tau)$ that underly the resulting analysis coincide with those for the single transformation $\tau^{t_0}$ for a suitable choice of time-step $t_0 > 0$).  As a result Potts is able to obtain convergence and multiple recurrence for the above continuous averages, as well as various finer results describing their limit.

A different approach to the continuous-time setting has also recently been investigated by Bergelson, Leibman and Moreira~\cite{BerLeiMor10}.  Their work is based on a direct reduction to analogous questions for discrete transformations, for which the desired convergence is then assumed as a black-box result, rather than showing how discrete-action machinery can be reconstructed in the continuous setting. This reduction, in turn, is based on a very general result asserting that convergence to a common limit for a family of discrete C\'esaro averages implies the similar convergence for an enveloping continuous C\'esaro average (Theorem 0.1 in their preprint). Their method gives results for the above one-parameter linear multiple averages, and also for polynomial multiple averages such as
\[\barint_0^T \prod_{i=1}^k (f_i\circ \tau^{p_i(t)})\,\d t\quad\hbox{for polynomials}\ p_i:\bbR\to\bbR.\]

However, neither of the above approaches is currently sufficient to prove convergence for general polynomial averages in an action of $\bbR^D$ for $D \geq 2$, because in both cases the authors make heavy appeal to results for the discrete-time setting that are already known: either by showing how structural results from that setting have analogs for continuous-time flows, or by simply reducing the desired convergence assertion to known facts about discrete transformations.  By contrast, one of the surprising features of Theorem~\ref{thm:main} is that the analysis of continuous-time polynomial multiple averages in higher-dimensional actions is actually much \emph{simpler} than its discrete-time relative, so that we are able to prove the theorem above while the behaviour of the analogous discrete-time averages remains largely mysterious.  Our proof uses the machinery of sated extensions, previously employed in~\cite{Aus--lindeppleasant1,Aus--lindeppleasant2} and related to earlier arguments from~\cite{Aus--nonconv} and~\cite{Hos09}, but which has not yet provided enough insight to handle the general case of Conjecture~\ref{conj}.  (We should also remark that another consequence of using extensions in proofs of convergence is that relatively little can then be deduced about the structure of the limit, whereas Potts and Bergelson, Leibman and Moreira do obtain some such descriptive results as well.)  The greater simplicity in the continuous-time setting results from the greater `smoothness' of the continuous averages alluded to above, which will be exploited concretely through the use of the fractional-power time-changes that were mentioned there.

\subsubsection*{Acknowledgements}

I am grateful to Vitaly Bergelson for several helpful suggestions, and to Bryna Kra for motivating me to return to this project after a period of neglect.  This research was supported by a Fellowship from Microsoft Corporation and a Fellowship from the Clay Mathematics Institute.

\section{Formulating the right question}\label{sec:right-ques}

It turns out that one of the chief obstacles to proving
Theorem~\ref{thm:main} is finding the right conjecture: the
assertion of Theorem~\ref{thm:main} is actually too weak. Instead we
will consider functions $\bphi_i$ (rather than polynomials
$\bf{p}_i$) that are vector-valued sums of maps with powers between
$0$ and $1$.  This will be made simpler by also allowing some greater flexibility in the location of the intervals over which we take our ergodic averages.

\begin{dfn}[Tempered sequence; tempered-uniform convergence]
A sequence of bounded open intervals $I_n \subset \bbR$
is \textbf{tempered} if $|I_n|\to\infty$ and there is some fixed $K
\geq 0$ such that $\rm{dist}(0,I_n) \leq K\cdot |I_n|$ for all
$n\geq 1$.

Given now a locally integrable map $v:[0,\infty)\to \frH$ into some
Hilbert space, the averages of $v$ over bounded intervals in $\bbR$ converge
\textbf{tempered-uniformly} if the sequence of averages
$\frac{1}{|I_n|}\int_{I_n}v(t)\,\d t$ converges as $n\to\infty$ for
any tempered sequence of intervals $I_n \subset [0,\infty)$.
\end{dfn}

In fact, tempered-uniform convergence is not really a new property of such a map $v$. Indeed, one has
\[\barint_{a_n}^{b_n}v(t)\,\d t = \frac{b_n}{b_n - a_n}\barint_0^{b_n}v(t)\,\d t - \frac{a_n}{b_n - a_n}\barint_0^{a_n}v(t)\,\d t\]
for any real $a_n < b_n$, and if the sequence of intervals $(a_n,b_n)$ is tempered then the coefficients $b_n/(b_n - a_n)$ and $a_n/(b_n - a_n)$ that appear here remain bounded with difference equal to $1$, so the tempered-uniform convergence of the averages $\barint_I v(t)\,\d t$ is actually equivalent to the Ces\`aro convergence of $v$.  (This also clearly implies that the limit is the same along any tempered-uniform sequence of intervals.)  The usefulness of the above definition is in providing a convenient handle for this convergence along sequences of intervals that may not be pinned to the origin.

\begin{lem}\label{lem:change-vars}
Suppose that $v:[0,\infty) \to \frH$ is a norm-continuous and
bounded function taking values in some Hilbert space. Then the
tempered-uniform convergence of the averages
\[\barint_I v(t)\,\d t\]
is equivalent to that of the averages
\[\barint_I v(s^\a)\,\d s\]
for any $\a > 0$, and the limits are the same.
\end{lem}

\textbf{Proof}\quad Let $I = (a,b)$.  By symmetry it suffices to
prove the forward implication and show that the limits are the same.
This is achieved by making the substitution $t := s^\a$ in the
second integral to obtain
\[\barint_a^b v(s^\a)\,\d s = \frac{1}{b-a}\int_{a^\a}^{b^\a} v(t) \frac{1}{\a}t^{1/\a - 1}\,\d t.\]

We will now show how this can always be written as a weighted
integral of averages of the form $A_J(v) := \barint_J v(t)\,\d t$ in
such a way that we can apply the tempered-uniform convergence of
these latter.  Let us write $A_\infty(v)$ for their tempered uniform
limit.  Since the case $\a = 1$ is trivial, this argument falls into
two remaining cases.

\quad\textbf{Case 1:}\quad$\underline{0 < \a < 1}$\quad In fact let
us assume further that $\a \neq \frac{1}{2}$, since this irritating
case can be treated by exactly the same method except that the
antiderivative of the function $t\mapsto t^{1/\a - 2}$ cannot be
written using the usual formula that is valid for other $\a$.

We first compute that
\begin{eqnarray*}
\frac{1}{b-a}\int_{a^\a}^{b^\a} v(t) \frac{1}{\a}t^{1/\a - 1}\,\d t
&=& \frac{\frac{1}{\a}(a^\a)^{1/\a - 1}}{b-a}\int_{a^\a}^{b^\a}
v(t)\,\d t\\
&& + \frac{1}{b-a}\int_{a^\a}^{b^\a} v(t)
\frac{1}{\a}\big(t^{1/\a - 1} - (a^\a)^{1/\a - 1}\big)\,\d t\\
&=& \frac{a^{1 - \a}(b^\a - a^\a)}{\a(b-a)}A_{(a^\a,b^\a)}(v)\\
&& + \frac{1}{b-a}\int_{a^\a}^{b^\a} \frac{1/\a - 1}{\a}t^{1/\a -
2}\int_t^{b^\a}v(u)\,\d u\,\d t\\
&=& \frac{a^{1 - \a}(b^\a - a^\a)}{\a(b-a)}A_{(a^\a,b^\a)}(v)\\
&& + \frac{1 - \a}{\a^2(b-a)}\int_{a^\a}^{b^\a} (b^\a - t)t^{1/\a -
2}A_{(t,b^\a)}(v)\,\d t.
\end{eqnarray*}

Now suppose that $I_n = (a_n,b_n)\subset [0,\infty)$ is a tempered
sequence of intervals, say with $a_n \leq K(b_n - a_n)$ for all $n$.
Define $K_n \in (0,1]$ by $b_n - a_n =: K_n b_n$, so that a re-arrangement of the temperedness inequality gives $K_n \geq \frac{1}{1 + K}$, and now a simple computation gives
\[\frac{\a}{1 - \a}(1 - (1 - K_n)^{1 - \a}) - \a K_n = \int_{(1 - K_n)^\a}^1 (1 - t)t^{1/\a - 2}\,\d t,\]
which is uniformly positive over the possible values $K_n \in [1/(1 + K),1]$.  Therefore for any $\eps > 0$ we can select some $L > K$ (typically much larger than $K$) such that
\begin{multline*}
\frac{\a}{1 - \a}\Big(\Big(\frac{L}{1 + L}\Big)^{1 - \a} - (1 - K_n)^{1 - \a}\Big) - \a \Big(K_n - \frac{1}{1 + L}\Big)\\
\geq (1 - \eps)\Big(\frac{\a}{1 - \a}(1 - (1 - K_n)^{1 - \a}) - \a K_n\Big)
\end{multline*}
for any $K_n \in [1/(1 + K),1]$.  Defining $c_n := (1 - \frac{1}{1 + L})b_n$, we now have $a_n < c_n < b_n$, but on the other hand $c_n = L(b_n - c_n)$, so any sequence of intervals of the form $(c'_n,b_n)$ or $((c'_n)^\a,b_n^\a)$ for some selection of $c_n' \in [a_n,c_n]$ is still tempered.

Using the above inequality, another simple calculation gives
\begin{eqnarray*}
&&\int_{a_n^\a}^{c_n^\a} (b_n^\a - t)t^{1/\a - 2}\,\d t\\ &&=
\frac{b_n^\a}{(1/\a - 1)}((c_n^\a)^{1/\a - 1} - (a_n^\a)^{1/\a - 1})
\\ &&\quad - \frac{1}{1/\a}((c_n^\a)^{1/\a} - (a_n^\a)^{1/\a})\\&&=
\frac{b_n^\a}{(1/\a - 1)}(c_n^{1 - \a} - a_n^{1 - \a}) - \frac{1}{1/\a}(c_n - a_n)\\
&&= b_n\Big(\frac{\a}{1 - \a}\Big(\Big(\frac{L}{1 + L}\Big)^{1 - \a} - (1 - K_n)^{1 - \a}\Big) - \a \Big(K_n - \frac{1}{1 + L}\Big)\Big)\\
&&\geq
(1 - \eps)b_n\Big(\frac{\a}{1 - \a}(1 - (1 - K_n)^{1 - \a}) - \a K_n\Big)\\
&&=(1 - \eps)\Big(\frac{b_n^\a}{(1/\a - 1)}(b_n^{1 - \a} - a_n^{1 - \a}) - \frac{1}{1/\a}(b_n - a_n)\Big)\\
&&= (1-\eps) \int_{a_n^\a}^{b_n^\a} (b_n^\a - t)t^{1/\a - 2}\,\d
t.
\end{eqnarray*}
Combining this with the change-of-variables made above and using
that $v$ is a bounded function we can write
\begin{multline*}
\barint_{a_n}^{b_n} v(s^\a)\,\d s = \frac{a_n^{1 - \a}(b_n^\a -
a_n^\a)}{\a(b_n-a_n)}A_{(a_n^\a,b_n^\a)}(v)\\ + \frac{1 -
\a}{\a^2(b_n-a_n)}\int_{a_n^\a}^{c_n^\a} (b_n^\a - t)t^{1/\a -
2}A_{(t,b_n^\a)}(v)\,\d t + R_1
\end{multline*}
where $\|R_1\|_2 \leq \eps\|v\|_{L^\infty([0,\infty)\to\frH)}$. Now
as $n\to\infty$, any sequence of intervals $I'_n$ such that $I'_n =
(t,b_n^\a)$ for some $a_n^\a \leq t \leq c_n^\a$ for each $n$ also
satisfies $\rm{dist}(0,I'_n) \leq K'|I'_n|$ for some $K' \geq 0$.
Hence the tempered-uniform convergence of the averages $A_I(v)$ to
some limiting vector, say $A_\infty(v)$, implies that we must
actually have
\[\|A_{(t,b_n^\a)}(v) - A_\infty(v)\| \leq \eps\|v\|_{L^\infty([0,\infty)\to\frH)}\quad\quad\forall t \in [a_n^\a,c_n^\a)\]
for all sufficiently large $n$.  Inserting this approximation into
the above average, we deduce that for all sufficiently large $n$ we
have
\begin{eqnarray*}
\barint_{a_n}^{b_n} v(s^\a)\,\d s &=& \frac{a_n^{1 - \a}(b_n^\a -
a_n^\a)}{\a(b_n-a_n)}A_\infty(v)\\ && \quad + \frac{1
- \a}{\a^2(b_n-a_n)}\int_{a_n^\a}^{c_n^\a} (b_n^\a - t)t^{1/\a -
2}A_\infty(v)\,\d t + R_2 + R_1\\
&=& \frac{a_n^{1 - \a}(b_n^\a - a_n^\a)}{\a(b_n-a_n)}A_\infty(v)\\
&& \quad + \frac{1 -
\a}{\a^2(b_n-a_n)}\int_{a_n^\a}^{b_n^\a} (b_n^\a - t)t^{1/\a -
2}A_\infty(v)\,\d t + R_3 + R_2 + R_1\\
&=& A_\infty(v) + R_3 + R_2 + R_1
\end{eqnarray*}
with $\|R_2\|,\|R_3\| \leq \eps\|v\|_{L^\infty([0,\infty)\to\frH)}$.
Since $\eps$ was arbitrary we conclude that
\[\barint_{a_n}^{b_n} v(s^\a)\,\d s\to A_\infty(v)\quad\quad\hbox{as}\ n\to\infty,\]
as required.

\quad\textbf{Case 2:}\quad$\underline{\a > 1}$\quad This is similar,
except now we start with the computation
\begin{eqnarray*}
\frac{1}{b-a}\int_{a^\a}^{b^\a} v(t) \frac{1}{\a}t^{1/\a - 1}\,\d t
&=& \frac{\frac{1}{\a}(b^\a)^{1/\a - 1}}{b-a}\int_{a^\a}^{b^\a}
v(t)\,\d t\\
&& + \frac{1}{b-a}\int_{a^\a}^{b^\a} v(t)
\frac{1}{\a}\big(t^{1/\a - 1} - (b^\a)^{1/\a - 1}\big)\,\d t\\
&=& \frac{b^{1 - \a}(b^\a - a^\a)}{\a(b-a)}A_{(a^\a,b^\a)}(v)\\
&& + \frac{1}{b-a}\int_{a^\a}^{b^\a} \frac{1 - 1/\a}{\a}t^{1/\a - 2}\int_{a^\a}^tv(u)\,\d u\,\d t\\
&=& \frac{b^{1 - \a}(b^\a - a^\a)}{\a(b-a)}A_{(a^\a,b^\a)}(v)\\
&& + \frac{\a - 1}{\a^2(b-a)}\int_{a^\a}^{b^\a} (t - a^\a)t^{1/\a -
2}A_{(a^\a,t)}(v)\,\d t.
\end{eqnarray*}
Just as before, given a tempered sequence of intervals $(a_n,b_n)$
we can approximate the above expressions arbitrarily well by
truncating the second integrals so that they involve only
$A_{(a_n^\a,t)}$ for $c_n^\a \leq t\leq b_n^\a$ for some suitable
sequence $a_n < c_n < b_n$, and then argue that all of the
expressions $A_{(a_n^\a,t)}$ for such $t$ converge uniformly fast to
$A_\infty(v)$, so that the left-hand integrals must do the same.
Note that in this case, if the values $a_n$ actually tend to $0$,
then this argument also requires the fact that $1/\a - 1 > -1$ and
so the function $t\mapsto t^{1/\a - 1}$ is locally integrable at
$0$: this is needed so that the small range $(a_n^\a,c_n^\a]$ that
we initially omit from the second integral on the right-hand side
can be chosen so as to give an arbitrarily small contribution
overall. \qed

Now let us examine a little the taxonomy of the functions that result from applying such fractional-power time-changes to polynomials.

\begin{dfn}[Fractional polynomial; height; degree; goodness]
If $d \in \bbN$ then a map $\bphi:[0,\infty)\to \bbR^D$ is a \textbf{fractional
polynomial} (`\textbf{f-polynomial}') of \textbf{height $d$} if it takes the form
\[\bphi(t) = \sum_{j=1}^d t^{j/d}\bf{v}_j\]
for some tuple of vectors $\bf{v}_1$, $\bf{v}_2$, \ldots,
$\bf{v}_d \in \bbR^D$. Note that there may be more than one possible choice of the height $d$ for a given map $\bphi$, and in the following we will need to keep track of the height as well.  We allow the possibility that $\bf{v}_d = \bs{0}$, and let the \textbf{degree} of $\bphi$, denoted by $\deg\bphi$, be the largest fraction $j/d$ for which $\bf{v}_j \neq \bs{0}$; we also say that $\bphi$ is of \textbf{top-degree} if $\deg\bphi = d/d = 1$.

A height-$d$ f-polynomial $\bphi$ is \textbf{good} if the list of vectors $\bf{v}_1$, $\bf{v}_2$, \ldots,
$\bf{v}_{d\cdot\deg\bphi}$ is linearly independent (including the assertion that they are all non-zero).  Note that this property depends on the choice of height: for example, given $\bphi$ as above, we may always re-write it as a sum $\sum_{j=1}^{2d} t^{j/2d}\bf{v}'_j$ in which $\bf{v}'_j = \bs{0}$ for all odd $j$.

Also, we set
\[V(\bphi) := \rm{span}\{\bf{v}_1,\bf{v}_2,\ldots,\bf{v}_d\},\]
which agrees with the linear span of the image $\phi(\bbR) \subseteq \bbR^d$, and we set
\[\bphi^-(t) := \sum_{j=1}^{d-1} t^{j/d}\bf{v}_j\]
(so $\bphi^- = \bphi$
when $\bphi$ is not of top-degree).
\end{dfn}

\begin{dfn}[Good families of f-polynomials]
A family $\vec{\bphi} = \{\bphi_1,\bphi_2,\ldots,\bphi_k\}$ of f-polynomials, say expressed as
\[\bphi_i = \sum_{j=1}^d t^{j/d}\bf{v}_{i,j}\quad\quad\hbox{for}\ i=1,2,\ldots,d,\]
is \textbf{good} if each of the individual f-polynomials $\bphi_i$ is good, and moreover all of the non-zero vectors $\bf{v}_{i,j}$ appearing in the above expressions are linearly independent.
\end{dfn}

We can now give the simple reformulation of Theorem~\ref{thm:main} in terms of fractional polynomials:

\begin{thm}\label{thm:A}
For any action $\tau:\bbR^D \curvearrowright (X,\mu)$, any good family of f-polynomials
$\bphi_1$, $\bphi_2$, \ldots, $\bphi_k:[0,\infty)\to\bbR^D$ and any tempered sequence of intervals $I_n \subset
[0,\infty)$ the associated averages
\[A_n(f_1,f_2,\ldots,f_k) := \barint_{I_n} \prod_{i=1}^k(f_i\circ \tau^{\bphi_i(t)})\,\d t\]
converge in $L^2(\mu)$ as $n\to\infty$ for any $f_1$, $f_2$, \ldots,
$f_k \in L^\infty(\mu)$.
\end{thm}

\textbf{Proof of Theorem~\ref{thm:main} from
Theorem~\ref{thm:A}}\quad This relies on a simple changed of
variables.  Suppose that $T:\bbR^D\actson (X,\mu)$, that
\[\bf{p}_i(t) = \sum_{j=1}^dt^j\bf{u}_{i,j}\quad\quad\hbox{for}\ i=1,2,\ldots,k,\]
and let $v(t) := \prod_{i=1}^k (f_i\circ \tau^{\bf{p}_i(t)})$, so this is a norm-continuous and bounded map into $L^2(\mu)$.

We first reduce the case in which all the $\bf{u}_{i,j}$ are non-zero and linearly independent. Introduce a formal collection of vectors $\bf{v}_{i,j}$ for $i\leq k$ and $j \leq d$ that are all non-zero and linearly independent.  Let $\bbR^{D'}$ be their linear span, and define a linear map $A:\bbR^{D'}\to \bbR^D$ by setting
\[A\bf{v}_{i,j} = \bf{u}_{i,j}\]
(so there may be linear dependences among these images, and some of them may be zero).  Now define $T':\bbR^{D'}\actson (X,\mu)$ by $(T')^{\bf{v}} := T^{A\bf{v}}$ and also
\[\bf{p}'_i(t) = \sum_{j=1}^dt^j\bf{v}_{i,j}\quad\quad\hbox{for each}\ i=1,2,\ldots,k,\]
and observe that $T^{\bf{p}_i(t)} = (T')^{\bf{p}_i'(t)}$.  Hence we can write our multiple averages in terms of $T'$ and the $\bf{p}_i'$, and so it suffices to prove Theorem~\ref{thm:main} in case the coefficient vectors $\bf{u}_{i,j}$ are all non-zero and linearly independent.

Assuming this, define
\[\bphi_i(s) := \sum_{j=1}^d s^{j/d}\bf{u}_{i,j}\quad\hbox{for}\ j=1,2,\ldots,d,\]
and observe that this is a good family.  On the other hand, by Lemma~\ref{lem:change-vars}, the norm convergence of the averages
$\barint_0^T v(t)\,\d t$ follows from the tempered-uniform
convergence of the averages
\[\barint_I v(s^{1/d})\,\d s = \barint_I \prod_{i=1}^k (f_i\circ \tau^{\bphi_i(s)})\,\d s.\]
This completes the proof. \qed

\section{The full induction}

\subsection{Furstenberg self-joinings and partially characteristic factors}

Here we reformulate in our present setting some older machinery from
the study of multiple ergodic averages.

First, suppose that we have already established convergence for some
family $\{\bphi_1,\bphi_2,\ldots,\bphi_k\}$ of f-polynomials $[0,\infty)
\to\bbR^D$, that $\tau:\bbR^D\curvearrowright(X,\mu)$ is a system
and that $A_0$, $A_1$, $A_2$, \ldots, $A_k$ are Borel subsets of
$X$.  Then our assumption implies that the scalar averages
\[\barint_0^T \mu(A_0\cap \tau^{-\bphi_1(t)}(A_1)\cap\cdots\cap \tau^{-\bphi_k(t)}(A_k))\,\d t\]
converge as $T\to\infty$.  Denoting the limiting value by
$\mu^\rm{F}_{\bphi_1,\bphi_2,\ldots,\bphi_k}(A_0\times
A_1\times \cdots \times A_k)$, it is now easy to check that this may
be extended by multilinearity and continuity, and that this actually
defines a $(k+1)$-fold self-joining
$\mu^\rm{F}_{\bphi_1,\bphi_2,\ldots,\bphi_k}$ of the system
$(X,\mu,\tau)$ (which depends also on $\tau$, although we suppress this in our notation).  This construct has its origins in Furstenberg's
original work on multiple recurrence~\cite{Fur77}, and is referred
to as the \textbf{Furstenberg self-joining} associated to the family
$\{\bphi_1,\bphi_2,\ldots,\bphi_k\}$.  This will prove an important
tool for analyzing our averages, much as in the previous works
in~\cite{ConLes84,Zie07,Aus--nonconv,Aus--lindeppleasant1}.

Next, a factor
$\xi:(X,\mu,\tau)\to (Y,\nu,\s)$ is
\textbf{partially characteristic} for a given tuple of f-polynomials
$\bphi_1$, $\bphi_2$, \ldots, $\bphi_k$ if
\[\barint_{I_n} \prod_{i=1}^k(f_i\circ \tau^{\bphi_i(t)})\,\d t \sim \barint_{I_n}  \prod_{i=1}^{k-1}(f_i\circ \tau^{\bphi_i(t)})\cdot (\sfE_\mu(f_k\,|\,\xi)\circ \tau^{\bphi_k(t)})\,\d t\]
for any $f_1$, $f_2$, \ldots, $f_k \in L^\infty(\mu)$ and any
tempered sequence of intervals $I_n \subset [0,\infty)$, where we
write $F_n \sim G_n$ to denote that $\|F_n - G_n\|_2 \to 0$ as
$N\to\infty$.  This notion is based on a definition that first
appears in Furstenberg and Weiss' paper~\cite{FurWei96}, and which
has gone through a number of incarnations since.  Note that it involves operating only on the last function $f_k$ in the list.

As in most previous proofs of convergence for some family of
multiple ergodic averages, the heart of our induction will be
finding a partially characteristic factor that has some additional
structure allowing averages of interest to be re-written into a
simpler form.  Here we will also make crucial use of a more recent
twist on this strategy, in which an initially-given system must
first be extended (that is, expressed as a factor of some `larger'
system) before the desired factor can be shown to be partially
characteristic.  This approach has been developed
in~\cite{Aus--nonconv,Aus--newmultiSzem,Aus--lindeppleasant1,Aus--lindeppleasant2}.
Here we will appeal to the very general machinery of sated
extensions from~\cite{Aus--lindeppleasant1} in order to make this
initial construction of an enlarged system.

\subsection{Satedness}

Following~\cite{Aus--lindeppleasant1}, a class $\sfC$ of standard
Borel probability-preserving $\bbR^D$-systems is sated to be
\textbf{idempotent} if it is closed under isomorphisms, inverse
limits and (not necessarily ergodic) joinings.  These conditions are
enough to guarantee that any such system $(X,\mu,\tau)$ has a
maximal factor $\zeta^\tau_\sfC:(X,\mu,\tau)\to \sfC(X,\mu,\tau)$
whose target system is a member of the class $\sfC$.  In these
terms, a system is \textbf{$\sfC$-sated} if whenever
$\pi:(\t{X},\t{\mu},\t{\tau}) \to (X,\mu,\tau)$ is an extension, the
factors $\pi$ and $\zeta^{\t{\tau}}_\sfC$ are relatively independent
over the further common factor $\zeta^{\tau}_\sfC\circ \pi$.  More concretely, the $\sfC$-satedness of $(X,\mu,\tau)$ means that if $f \in L^2(\mu)$ and we prove that
\[\sfE_\mu(f\circ \pi\,|\,\zeta^{\t{\tau}}_\sfC) \neq 0\]
for some extension $\pi:(\t{X},\t{\mu},\t{\tau})\to (X,\mu,\tau)$, then we can deduce that in fact
\[\sfE_\mu(f\,|\,\zeta^\tau_\sfC)\neq 0.\]

Here we will be concerned with idempotent classes of the following
kinds. If $V\leq \bbR^D$ is a vector subspace then we may associate
to it the class $\sfZ_0^V$ of $\bbR^D$-systems whose $V$-subaction
is trivial; and now given several subspaces $V_1,V_2,\ldots,V_\ell
\leq \bbR^D$ we let $\sfZ_0^{V_1}\vee \sfZ_0^{V_2}\vee \cdots \vee
\sfZ_0^{V_\ell}$ be the class of all joinings of systems drawn from
each of the classes $\sfZ_0^{V_i}$.  Both of these examples are
readily verified to be idempotent (or see Section 3 in~\cite{Aus--lindeppleasant1}), and if we abbreviate
$\zeta_{\sfZ_0^V}^\tau =: \zeta_V^\tau$ then a simple check also
shows that
\[\zeta_{\sfZ_0^{V_1}\vee\sfZ_0^{V_2}\vee\cdots\vee\sfZ_0^{V_\ell}}^\tau \simeq \bigvee_{i=1}^\ell\zeta_{V_i}^\tau\]
as factor maps of $(X,\mu,\tau)$, in the sense that these maps define the same factor of $\tau$ up to negligible sets, and where the right-hand side denotes the factor map generated by the individual factor maps $\zeta^\tau_{V_i}$.

In these terms, we now make the following definition, which closely
follows the idea of a `fully isotropy-sated' system
from~\cite{Aus--lindeppleasant1}.

\begin{dfn}[Fully rationally sated systems]\label{dfn:FRS}
We will write that $(X,\mu,\tau)$ is \textbf{fully rationally sated}
(`\textbf{FRS}') if it is sated for the idempotent classes
\[\bigvee_{V \in \cal{V}}\sfZ_0^V\]
whenever $\cal{V}$ is a finite collection of rational subspaces of $\bbR^D$, where a subspace $V\leq \bbR^D$ is \textbf{rational} if it has a
basis consisting of members of $\bbQ^D$.
\end{dfn}

The usefulness of this definition derives from the following general
fact, which is an immediate special case of Theorem 3.11
in~\cite{Aus--lindeppleasant1}:

\begin{thm}[Sated extensions exist]\label{thm:FRS}
Every $\bbR^d$-system has an extension that is FRS. \qed
\end{thm}

The principal limitation of Theorem 3.11 in~\cite{Aus--lindeppleasant1} is that it can be applied to at most countably many families of idempotent classes, and this is why we must
restrict our attention to rational subspaces in
Definition~\ref{dfn:FRS}.

\subsection{An ordering on fractional polynomial families}\label{subs:prec}

The proof of Theorem~\ref{thm:A} will require an induction on good families of f-polynomials $\vec{\bphi}$, and so we make a separate step of introducing the relevant ordering on these families, and introducing two particular kinds of `downward movement' through the collection of such families that will appear during the induction.

\begin{dfn}[The precedence ordering]
Given two non-empty good families of height-$d$ fractional polynomials, say $\vec{\bphi} = \{\bphi_1,\bphi_2,\ldots,\bphi_k\}$ and $\vec{\bpsi} = \{\bpsi_1,\bpsi_2,\ldots,\bpsi_\ell\}$, we say that $\vec{\bpsi}$ \textbf{precedes} $\vec{\bphi}$, written
$\vec{\bpsi} \prec \vec{\bphi}$, if
\begin{itemize}
\item $\ell \leq k$,
\item when the $\bphi_i$ and $\bpsi_i$ are ordered so that their degrees are non-increasing in $i$, one has
\[\deg\bpsi_i \leq \deg\bphi_i\quad\forall i \leq \ell\]
(so the $\bphi_i$ for $\ell < i \leq k$ are not needed here), and
\item either $\ell < k$, or if $\ell = k$ then strict inequality holds in the second condition above for some $i\leq \ell$.
\end{itemize}

In addition, we always have $\vec{\bpsi} \prec \vec{\bphi}$ if $\vec{\bpsi}$ and $\vec{\bphi}$ are families of fractional polynomials of distinct heights $d$ and $d' > d$.
\end{dfn}

It is clear that this defines a partial order on the collection of all families of height-$d$ fractional polynomials, and that it satisfies the descending chain condition.  Note that the inequality $\ell < k$ is not by itself enough to guarantee that $\vec{\bpsi}\prec \vec{\bphi}$: it is also necessary that the former family not have too many high-degree members, in the sense of the second point above.

Suppose now that $\vec{\bphi} = \{\bphi_1,\bphi_2,\ldots,\bphi_k\}$ is a good family with $k\geq 2$.  Two special kinds of precedent for $\vec{\bphi}$ will be important in the sequel:
\begin{itemize}
\item On the one hand, suppose that $j_1 \in \{1,2,\ldots,d\}$ is minimal such that there is some $\bphi_{i_1} \in \vec{\bphi}$ for which the leading term has degree $j_1/d$. In this case the family
\[\vec{\bpsi} := \{\bphi_i - \bphi_{i_1}:\ i \in \{1,2,\ldots,k\}\setminus \{i_1\}\}\]
precedes $\vec{\bphi}$.  Indeed, we have removed one instance of degree $j_1/d$, and for every $i \neq i_1$ we have that $\bphi_i - \bphi_{i_1}$ still has degree at least $\deg\bphi_i \geq j_1/d$ (the goodness assumption implies that all coefficients of $\bphi_i$ and $\bphi_{i_1}$ are distinct, so there can be no cancellation of coefficients here, and $\deg\bphi_i \geq j_1/d$ by the minimality of $j_1$).  Also, $\vec{\bpsi}$ is still a good family since the collection of its coefficients is
\[\{\bf{v}_{i,j}:\ i \in \{1,2,\ldots,k\}\setminus\{i_1\}\}\quad\hbox{for}\ j \in \{j_1 + 1,\ldots,d\}\]
and
\[\{\bf{v}_{i,j} - \bf{v}_{i_1,j}:\ i \in \{1,2,\ldots,k\}\setminus\{i_1\}\}\quad\hbox{for}\ j \leq j_1,\]
and it is clear that these are all still nonzero and linearly independent.

We refer to a family $\vec{\bpsi}$ constructed this way as a \textbf{precedent of $\vec{\bphi}$ of type I}.

\item On the other hand, if $\bphi_i\in \vec{\bphi}$ is a top-degree member then the family
\[\{\bphi_1,\bphi_2,\ldots,\bphi_{i-1},\bphi_i^-,\bphi_{i+1},\ldots,\bphi_k\}\]
has swapped out an entry of degree $\deg\bphi_i$ and replaced it with an entry of degree $\deg\bphi_i - 1$ (where in case $d = 1$ we instead take the above to mean that $\bphi_i$ has been omitted altogether). So again it clearly precedes $\vec{\bphi}$ and is still good.  We call this a \textbf{precedent of $\vec{\bphi}$ of type II}.
\end{itemize}

\subsection{The main induction}

In order to formulate an inductive hypothesis that includes Theorem~\ref{thm:A} and can be closed on itself, we will actually prove a composite of three different
properties of the multiple averages associated to each family
$\vec{\bphi}$.  To this end we insert the statement of
Theorem~\ref{thm:A} as the second of the three related conclusions.

\begin{thm}\label{thm:B}
Suppose that $(X,\mu,\tau)$ is an $\bbR^D$-system and that $\vec{\bphi} = \{\bphi_1, \bphi_2, \ldots, \bphi_k\}$ is a good family of height-$d$ f-polynomials $\bbR\to \bbR^D$ with expressions
\[\bphi_i(t) = \sum_{j=1}^dt^{j/d}\bf{v}_{i,j},\]
where each $\bf{v}_{i,j} \in \bbQ^D$. Then the following hold:
\begin{itemize}
\item[A($\vec{\bphi}$):] if $(X,\mu,\tau)$ is FRS and $\bphi_k \in \vec{\bphi}$ is of top-degree then the factor
\[\xi := \zeta^\tau_{\bbR\bf{v}_{k,d}}\vee\bigvee_{i = 1}^{k-1}\zeta^\tau_{V(\bphi_i - \bphi_k)}\]
is partially characteristic for the averages
\[A_I(f_1,f_2,\ldots,f_k) := \barint_I \prod_{i=1}^k f_i\circ \tau^{\bf{p}_i(t)}\,\d t\]
associated to $\vec{\bphi}$;
\item[B($\vec{\bphi}$):] the averages $A_I(f_1,f_2,\ldots,f_k)$ converge tempered-uniformly in
$L^2(\mu)$ for any $f_1,f_2,\ldots,f_k \in L^\infty(\mu)$;
\item[C($\vec{\bphi}$):] the Furstenberg self-joining
$\mu^\rmF_{\id,\bphi_1,\bphi_2,\ldots,\bphi_k}$ exists and is
invariant under the \textbf{off-diagonal flow}
\[t\mapsto \id\times \tau^{t\bf{v}_{1,j}}\times \tau^{t\bf{v}_{2,j}}\times \cdots\times \tau^{t\bf{v}_{k,j}}\]
for each $j\leq d$.
\end{itemize}
\end{thm}

Clearly the assumption that $\bf{v}_{i,j} \in \bbQ^D$ involves no loss of generality in that we may simply change basis to make it true, but it does fix what we mean by `FRS' in the statement of conclusion $\rm{A}(\vec{\bphi})$, and for this it is important.

In the above notation, for $\rm{D} \in \{\rm{A},\rm{B},\rm{C}\}$
we will write $\rm{D}(\prec\!\vec{\bphi})$ for the assertion that
$\rm{D}(\vec{\bpsi})$ holds for all $\vec{\bpsi} \prec
\vec{\bphi}$.

Our next task is to establish the base case of Theorem~\ref{thm:B}.

\begin{lem}~\label{lem:base-case}
The conclusion $\rm{A}(\vec{\bphi})\vee \rm{B}(\vec{\bphi})\vee \rm{C}(\vec{\bphi})$ holds if $\vec{\bphi}$ contains only one member, say $\bphi_1$.
\end{lem}

\textbf{Proof}\quad In this case we have
\[A_I = \barint_I f\circ \tau^{\bphi_1(t)}\,\d t\]
with $\bphi_1(t) = \sum_{j=1}^d t^{j/d}\bf{v}_j$.  Clearly by omitting extraneous directions we may assume that $\bbR^D = \bbR^d = V(\bphi_1)$.  Applying Lemma~\ref{lem:change-vars} in reverse, we see that it suffices instead to work with the polynomial $\bf{p}_1(t) = \sum_{j=1}^d t^j\bf{v}_j$.  However, for this we can simply regard the action $\tau$ as a unitary flow $U$ on the Hilbert space $L^2(\mu)$ using the Koopman representation, so that the above average (with $\bf{p}_1$ in place of $\bphi_1$) becomes
\[\barint_I U^{\bf{p}_i(s)}f\,\d s.\]
We can now show that this converges in norm to the conditional expectation $\sfE_\mu(f\,|\,\zeta_{\bbR^D}^\tau)$, which agrees with the orthogonal projection of $f$ onto the subspace of $U$-invariant functions in $L^2(\mu)$.  This is most easily proved by appeal to the Spectral Theorem, which translates this assertion into the fact that for any finite measure $\theta$ on $\widehat{\bbR^D}$ (which arises as the spectral measure assocated to $f$) the functional averages
\[\chi \mapsto \barint_I \chi(\bf{p}_1(s))\,\d s\]
converge tempered-uniformly to the indicator function at the origin, $1_{\{\bs{0}\}}$, in $L^2(\theta)$.  This in turn is clear using the Dominated Convergence Theorem, because these functions are all bounded by $1$ and if $\chi = (\chi_1,\chi_2,\ldots,\chi_d) \neq \bs{0}$ then the above average is
\[\barint_I \prod_{j=1}^d \chi_j(s^j)\,\d s,\]
whose convergence to $0$ uniformly as $|I|\to \infty$ (without any assumption of temperedness) follows by Weyl's Equidistribution Theorem (see, for example, Kuipers and Niederreiter~\cite{KuiNie74}).

This already contains the conclusion $\rm{B}(\vec{\bphi})$, and it also contains $\rm{A}(\vec{\bphi})$ because $\zeta^\tau_{\bbR^D}$ is contained in $\xi$ as a factor map.  Finally, $\rm{C}(\vec{\bphi})$ follows by applying this convergence inside the definition
\begin{multline*}
\mu^\rmF_{\id,\bphi_1}(A_0\times A_1) := \lim_{T\to\infty}\barint_0^T \mu(A_0\cap \tau^{-\bphi_1(t)}(A_1))\,\d t\\ = \int_X 1_{A_0}\cdot \Big(\lim_{T\to\infty}\barint_0^T 1_{A_1}\circ \tau^{\bphi_1(t)}\,\d t\Big)\,\d\mu,
\end{multline*}
since this now gives that $\mu^\rmF_{\id,\bphi_1}$ is the relatively independent product of two copies of $\mu$ over the factor $\zeta^\tau_{\bphi_1}$, which is invariant under everything. \qed

It remains to show the induction step.  Three different arguments will be used to complete one iteration of this, and we present these separately.

\begin{prop}\label{prop:C-A}
If $\vec{\bphi}$ contains at least two members then
\[\rm{C}(\prec\!\vec{\bphi}) \Longrightarrow \rm{A}(\vec{\bphi}).\]
\end{prop}

\textbf{Proof}\quad This is the most difficult of the three parts of our inductive implications. It will rest mainly on a reduction from $\vec{\bphi}$ to a precedent of type I.

Suppose that $f_1$, $f_2$, \ldots, $f_k
\in L^\infty(\mu)$ are such that
\[A_{I_n}(f_1,f_2,\ldots,f_k)\not\to 0\]
in $L^2(\mu)$ for some tempered sequence of intervals $I_n \subset
[0,\infty)$.  We will show that in this case we must also have
$\sfE_\mu(f_k\,|\,\xi) \neq 0$; from that point the proof is completed simply by substituting
the decomposition $f_k = \sfE_\mu(f_k\,|\,\xi) + (f_k -
\sfE_\mu(f_k\,|\,\xi))$ for a general function $f_k$ and appealing
to linearity in $f_k$.

Suppose that $\bphi_i(t) = \sum_{j=1}^d t^{j/d}\bf{v}_{i,j}$ for
$i=1,2,\ldots,k$ and that $j_1 \in \{1,2,\ldots,d\}$ is minimal such that some element of $\vec{\bphi}$ has degree $j_1/d$; by re-ordering we may assume that it is $\bphi_1$.

We begin by applying the version of the van der Corput estimate for
continuous families of vectors (see, for instance, Potts' Appendix B in~\cite{Pot09}) to the functions
\[\prod_{i=1}^kf_i\circ \tau^{\bphi_i(t)} \in L^2(\mu),\]
to deduce from the non-convergence
$A_{I_n}(f_1,f_2,\ldots,f_k)\not\to 0$ that we also have
\[\barint_0^H\barint_{I_n}\int_X \prod_{i=1}^k(f_i\circ \tau^{\bphi_i(t+h)})\cdot \prod_{i=1}^k(\overline{f_i}\circ \tau^{\bphi_i(t)})\ \d\mu\,\d t\,\d h \not\to 0\]
as $n\to\infty$ and then $H\to\infty$.

We may write the above integral as
\begin{multline*}
\barint_0^H\barint_{I_n}\int_X
\prod_{i=1}^k\big(\big((f_i\circ\tau^{\bphi_i(t+h) - \bphi_i(t)})
\cdot \overline{f_i}\big)\circ
\tau^{\bphi_i(t)}\big)\ \d\mu\,\d t\,\d h\\
=\barint_0^H\barint_{I_n}\int_X
\prod_{i=1}^k\big(\big((f_i\circ\tau^{\bphi^-_i(t+h) -
\bphi^-_i(t)}\tau^{h\bf{v}_{i,d}}) \cdot \overline{f_i}\big)\circ
\tau^{\bphi_i(t)}\big)\ \d\mu\,\d t\,\d h,
\end{multline*}
and now use the $\mu$-preserving transformation $\tau^{\bphi_1(t)}$
to change variables to obtain
\[\barint_0^H\barint_{I_n}\int_X
\prod_{i=1}^k\big(\big((f_i\circ\tau^{\bphi^-_i(t+h) -
\bphi^-_i(t)}\tau^{h\bf{v}_{i,d}}) \cdot \overline{f_i}\big)\circ
\tau^{\bphi_i(t) - \bphi_1(t)}\big)\ \d\mu\,\d t\,\d h.\] Now, since
$\bphi^-_i(t+h) - \bphi^-_i(t) \to 0$ for any $h \in \bbR$ as
$t\to\infty$, by the strong continuity of $\tau$ we also have
\[f_i\circ\tau^{\bphi_i^-(t+h) - \bphi^-_i(t)}\tau^{h\bf{v}_{i,d}} = f_i\circ \tau^{h\bf{v}_{i,d}} + \big(f_i\circ\tau^{\bphi_i^-(t+h) - \bphi^-_i(t)} - f_i\big)\circ \tau^{h\bf{v}_{i,d}}\to f_i\circ\tau^{h\bf{v}_{i,d}}\]
in $\|\cdot\|_2$ for any $h$ as $t\to\infty$ (and hence for $(1 -
\rm{o}(1))$-proportion of the interval $I_n$ as $n\to\infty$). Since all of these functions are also uniformly bounded in $\|\cdot\|_\infty$, replacing them one-by-one implies that as $n\to\infty$ the above integrals behave
asymptotically the same as
\[\barint_0^H\barint_{I_n}\int_X \prod_{i=1}^k\big(\big((f_i\circ \tau^{h\bf{v}_{i,d}})\cdot \overline{f_i}\big)\circ \tau^{\bphi_i(t) - \bphi_1(t)}\big)\ \d\mu\,\d t\,\d h.\]

Now consider the family $\vec{\bpsi} = \{\bpsi_2, \bpsi_3,
\ldots, \bpsi_k\}$ defined by
\[\bpsi_i(t) := \bphi_i(t) - \bphi_1(t),\]
which is a type-I precedent of $\vec{\bphi}$.  The inductive assumption of
$\rm{C}(\vec{\bpsi})$ implies that the above integrals converge
as $n\to\infty$ to
\[\barint_0^H\int_{X^k} \bigotimes_{i=1}^k\big((f_i\circ \tau^{h\bf{v}_{i,d}})\cdot \overline{f_i}\big)\ \d\mu^\rmF_{\id,\bpsi_2,\ldots,\bpsi_k}\,\d h,\]
and that the Furstenberg self-joining
$\mu^\rmF_{\id,\bpsi_2,\bpsi_2,\ldots,\bpsi_k} =:
\mu^\rmF_{\id,\vec{\bpsi}}$ that appears here is
invariant under each of the flows
\[t \mapsto \id\times \tau^{t(\bf{v}_{2,j} - \bf{v}_{1,j})}\times\tau^{t(\bf{v}_{3,j} - \bf{v}_{1,j})}\times\cdots\times \tau^{t(\bf{v}_{k,j} - \bf{v}_{1,j})}\]
for $j=1,2,\ldots,d$.

We now invoke the sleight-of-hand that underlies the previous works~\cite{Aus--nonconv,Aus--lindeppleasant1,Aus--lindeppleasant2} on nonconventional averages: we define a new jointly measurable action $\s:\bbR^D\actson X^k$ by setting
\[\s^{t\bf{v}_{i,j}} := \big(\tau^{t\bf{v}_{i,j}}\big)^{\times k}\quad\hbox{whenever}\ i\leq k-1\ \hbox{and}\ \bf{v}_{i,j}\neq \bs{0}\]
and
\[\s^{t\bf{v}_{k,j}} := \tau^{t\bf{v}_{1,j}}\times \tau^{t\bf{v}_{2,j}}\times \cdots\times \tau^{t\bf{v}_{k,j}}\quad\hbox{for all}\ j\leq d\]
(noting that $\bf{v}_{k,j}\neq \bs{0}$ for all $j\leq d$, by the assumptions that $\bphi_k$ is of top-degree and is good).  There is no conflict in this definition, because by the goodness assumption on $\vec{\bphi}$ all the non-zero vectors $\bf{v}_{i,j}$ are linearly independent.

In terms of $\s$, the measure $\mu^\rmF_{\id,\vec{\bpsi}}$ is invariant under $\s^{t\bf{v}_{i,j}}$ for any $i\leq k-1$ (since this is just an element of the diagonal action $\tau^{\times k}$), and also under each of the flows $\s^{t\bf{v}_{k,j}}$ for $j\leq d$, as a consequence of hypothesis $\rm{C}(\vec{\bpsi})$ and the identities
\[\s^{t\bf{v}_{k,j}} = \big(\id\times \tau^{t(\bf{v}_{2,j} - \bf{v}_{1,j})}\times\tau^{t(\bf{v}_{3,j} - \bf{v}_{1,j})}\times\cdots\times \tau^{t(\bf{v}_{k,j} - \bf{v}_{1,j})}\big)\circ \s^{t\bf{v}_{1,j}}\]

In particular, since $\mu^\rmF_{\id,\vec{\bpsi}}$ is $\s^{t\bf{v}_{k,d}}$-invariant, the usual Mean Ergodic Theorem implies that that the averages over $h \in (0,H)$ of the functions
\[\bigotimes_{i=1}^k (f_i\circ \tau^{h\bf{v}_{i,d}})  = \Big(\bigotimes_{i=1}^k f_i\Big)\circ \s^{h\bf{v}_{k,d}}\in L^\infty(\mu^\rmF_{\id,\vec{\bpsi}})\]
converge in $\|\cdot\|_2$ as $H\to\infty$ to some
$g \in L^\infty(\mu^\rmF_{\id,\vec{\bpsi}})$ that is
$\s^{\bbR\bf{v}_{k,d}}$-invariant.

Inserting this convergence into our integral averages, we deduce that they must tend to the limit
\[\int_{X^k} g\cdot\Big(\bigotimes_{i=1}^k\overline{f_i}\Big)\, \d\mu^\rmF_{\id,\vec{\bpsi}},\]
and therefore that this latter is not zero.

On the other hand, writing $\pi_i:X^k \to X$ for the $i^{\rm{th}}$
coordinate-projection, we see that for each $i = 1,2,\ldots,k-1$ and
each $j$ for which $\bf{v}_{i,j}\neq \bs{0}$, the function $f_i\circ \pi_i\in
L^\infty(\mu^\rmF_{\id,\vec{\bpsi}})$ is manifestly invariant under
the flow
\[t \mapsto \tau^{t(\bf{v}_{1,j} - \bf{v}_{i,j})}\times
\tau^{t(\bf{v}_{2,j} - \bf{v}_{i,j})}\times\cdots\times \tau^{t(\bf{v}_{k,j} -
\bf{v}_{i,j})} =
\s^{t(\bf{v}_{k,j}-\bf{v}_{i,j})},\] because the action of this flow on the $i^{\rm{th}}$ coordinate is trivial.

It follows that the $\bbR^D$-system
$(X^k,\mu^\rmF_{\id,\vec{\bpsi}},\s)$ is an extension of $(X,\mu,\tau)$ through the factor map $\pi_k$, and that among functions on this extended system we have that
$g$ is measurable with respect to
$\zeta^{\s}_{\bbR\bf{v}_{k,d}}$ and that $f_i\circ \pi_i$ is measurable with respect
to $\zeta^{\s}_{\bbR(\bf{v}_{i,j} - \bf{v}_{k,j})}$ for each
$i\leq k - 1$ and $j\leq d$, and hence with respect to $\zeta^{\s}_{V(\bphi_i -
\bphi_k)}$.  Therefore the non-vanishing of the above integral tells
us that the lifted function $f_k\circ \pi_k$ on the extended system
$(X^k,\mu^\rmF_{\id_X,\vec{\bpsi}},\s)$ has nonzero
conditional expectation onto the factor
\[\t{\xi} := \zeta^{\s}_{\bbR\bf{v}_{k,d}}\vee \bigvee_{i=1}^{k-1}\zeta^{\s}_{V(\bphi_i -
\bphi_k)}.\] Since $(X,\mu,\tau)$ is sated for the idempotent class
of all systems that can appear as such a factor (recall that one of the assumptions for assertion $\rm{A}(\vec{\bphi})$ is that $(X,\mu,\tau)$ is FRS), it follows that in
fact we must have $\sfE_\mu(f_k\,|\,\xi)\neq 0$, as required.  \qed

\begin{prop}\label{prop:AB-B}
If $\vec{\bphi}$ contains at least two members then
\[\rm{A}(\vec{\bphi})\vee\rm{B}(\prec\!\vec{\bphi}) \Longrightarrow
\rm{B}(\vec{\bphi}).\]
\end{prop}

\textbf{Proof}\quad First note that if $d'/d < 1$ is the maximal degree of any member of $\vec{\bphi}$, then performing the time-change $s := t^{d'/d}$ defines a new family of f-polynomials $\bphi_i'(s) := \bphi_i(s^{d/d'})$ of height $d'$ which has coefficient vectors precisely the same as for $\vec{\bphi}$, so it is still good.  By Lemma~\ref{lem:change-vars} it suffices to prove the desired convergence for the averages associated to these new f-polynomials $\bphi_i'$, but this is now an instance of $\rm{B}(\prec\!\vec{\bphi})$ because $d' < d$ and so this time-changed family precedes $\vec{\bphi}$.

Therefore we may assume that $\vec{\bphi}$ has some member of top-degree, say $\bphi_k$. In this case the result follows from Theorem~\ref{thm:FRS} using a reduction to a precedent of type II.

Next, it is clear that the averages $A_I$ associated to the family
$\vec{\bphi}$ converge for a system $(X,\mu,\tau)$ if they converge
for some extension $\pi:(\t{X},\t{\mu},\t{\tau})\to (X,\mu,\tau)$ of
that system, simply by lifting each function $f_i$ to $f_i\circ
\pi$, so in proving $\rm{B}(\vec{\bphi})$ we may assume that
$(X,\mu,\tau)$ is FRS.  Given this, $\rm{A}(\vec{\bphi})$
tells us that it suffices to prove convergence when $f_k$ is
measurable with respect to the factor
\[\xi := \zeta^\tau_{\bbR\bf{v}_{k,d}}\vee\bigvee_{i=1}^{k-1}\zeta^\tau_{V(\bphi_i - \bphi_k)}.\]
This, in turn, implies that it can be approximated in
$\|\cdot\|_2$ by a finite linear
combination of functions of the form $g\cdot \prod_{j=2}^k h_j$ with
$g$ invariant under $\tau^{\uhr \bbR\bf{v}_{k,d}}$ and each $h_i$
invariant under $\tau^{\uhr V(\bphi_i - \bphi_k)}$, where $\tau^{\uhr V}$ denotes the restriction of $\tau$ to an action of the subgroup $V$.  A simple
continuity argument now shows that it suffices to prove convergence
when $f_k$ equals a finite linear combination of such products, and
now by multilinearity it actually suffices to consider just
one such product.

However, having reduced to this case our averages enjoy the simplification
\begin{eqnarray*}
\barint_I \prod_{i=1}^k (f_i\circ \tau^{\bphi_i(t)})\,\d t &=& \barint_I \prod_{i=1}^{k-1} (f_i\circ \tau^{\bphi_i(t)})\cdot((gh_2h_3\cdots h_k)\circ \tau^{\bphi_k(t)})\,\d t\\
&=& \barint_I \prod_{i=1}^{k-1} ((h_if_i)\circ
\tau^{\bphi_i(t)})\cdot (g\circ \tau^{\bphi_k^-(t)})\,\d t,
\end{eqnarray*}
since $g\circ \tau^{\bphi_k(t)} =
g\circ \tau^{\bphi_k^-(t)}$ for all $t$, and $\bphi_i(t) - \bphi_k(t) \in V(\bphi_i -
\bphi_k)$ and so $h_i\circ \tau^{\bphi_k(t)} = h_i\circ
\tau^{\bphi_i(t)}$ for all $t$ and all $i \leq k-1$.

If $\bphi_k^-$ is nontrivial the convergence of these latter averages in $L^2(\mu)$ follows from that for a system of averages associated to the f-polynomial family $\{\bphi_1,\bphi_2,\ldots,\bphi_{k-1},\bphi_k^-\}$.  This is a type-II precedent of $\vec{\bphi}$, and so this convergence follows from the inductive assumption of $\rm{B}(\prec\vec{\bphi})$. On the other hand, if $\bphi_k^-$ is trivial (which can occur only if all of the maps $\bphi_i$ were actually linear), then we may simply move $g$ outside the above integral altogether and simply apply hypothesis $\rm{B}$ for the precedent family $\{\bphi_1,\bphi_2,\ldots,\bphi_{k-1}\}$. \qed

\begin{prop}\label{prop:ABC-C}
If $\vec{\bphi}$ has at least two members then
\[\rm{A}(\vec{\bphi})\vee\rm{B}(\vec{\bphi})\vee\rm{C}(\prec\!\vec{\bphi})
\Longrightarrow \rm{C}(\vec{\bphi}).\]
\end{prop}

\textbf{Proof}\quad Once again we may assume that $\vec{\bphi}$ has a member of top-degree, say $\bphi_k$, and in this case will use another reduction to a precedent of type II.

The Furstenberg self-joining
$\mu^\rm{F}_{\id,\vec{\bphi}}$ exists provided the averages
\[\barint_0^T \mu(A_0\cap \tau^{-\bphi(t)}(A_1)\cap \tau^{-\bphi_2(t)}(A_2)\cap\cdots\cap\tau^{-\bphi_k(t)}(A_k))\,\d t\]
converge as $T\to\infty$ for all Borel subsets $A_0$, $A_1$, \ldots,
$A_k \subseteq X$, and this now follows at once from
$\rm{B}(\vec{\bphi})$ by writing the above scalar-valued averages
as
\[\int_X 1_{A_0}\cdot A_{(0,T)}(1_{A_1},1_{A_2},\ldots,1_{A_k})\,\d\mu.\]

It remains to prove the invariance of
$\mu^\rmF_{\id,\vec{\bphi}}$ under the asserted off-diagonal
flows. If $\pi:(\t{X},\t{\mu},\t{\tau})\to (X,\mu,\tau)$ is any
extension then clearly $\mu^\rmF_{\id_X,\vec{\phi}}$ is the
pushforward of $\t{\mu}^\rmF_{\id_X,\vec{\phi}}$ under
$\pi^{\times (k+1)}$, and $\pi^{\times (k+1)}$ also intertwines each transformation
\[\id_X\times \t{\tau}^{t\bf{v}_{1,j}}\times \cdots\times \t{\tau}^{t\bf{v}_{k,j}}\quad\hbox{with}\quad \id_X\times \tau^{t\bf{v}_{1,j}}\times \cdots\times \tau^{t\bf{v}_{k,j}},\]
so by Theorem~\ref{thm:FRS} it will suffice to prove this invariance
assuming that $(X,\mu,\tau)$ is FRS.

Thus, suppose that $f_0$, $f_1$, $f_2$, \ldots, $f_k \in
L^\infty(\mu)$, $j\leq d$ and $t\in\bbR$. For each $j\leq d$ we will show that
\begin{multline*}
\int_{X^{k+1}}f_0\otimes f_1\otimes f_2\otimes \cdots\otimes
f_k\,\d\mu^\rm{F}_{\id,\vec{\bphi}}\\ = \int_{X^{k+1}}f_0\otimes
(f_1\circ \tau^{t\bf{v}_{1,j}})\otimes (f_2\circ
\tau^{t\bf{v}_{2,j}})\otimes \cdots\otimes (f_k\circ
\tau^{t\bf{v}_{k,j}})\,\d\mu^\rm{F}_{\id,\vec{\bphi}},
\end{multline*} from whence multilinearity and a simple approximation
argument complete the proof.

By the
hypothesis $\rm{A}(\vec{\bphi})$ we know that the factor
\[\xi := \zeta_{\bbR\bf{v}_{k,d}}^\tau\vee \bigvee_{i=1}^{k-1}\zeta^\tau_{V(\bphi_k - \bphi_i)}\]
is partially characteristic for the averages $A_{(0,T)}$ associated to
$\vec{\bphi}$. Hence the expression for the above integral as a
limit of integrals of $f_0\cdot A_{(0,T)}(f_1,f_2,\ldots,f_k)$ gives
that
\[\int_{X^{k+1}}f_0\otimes f_1\otimes f_2\otimes \cdots\otimes
f_k\,\d\mu^\rm{F}_{\id,\vec{\bphi}} = \int_{X^{k+1}}f_0\otimes
f_1\otimes f_2\otimes \cdots\otimes
\sfE_\mu(f_k\,|\,\xi)\,\d\mu^\rm{F}_{\id,\vec{\bphi}}\] and
similarly
\begin{multline*}
\int_{X^{k+1}}f_0\otimes (f_1\circ \tau^{t\bf{v}_{1,j}})\otimes
(f_2\circ \tau^{t\bf{v}_{2,j}})\otimes \cdots\otimes (f_k\circ
\tau^{t\bf{v}_{k,j}})\,\d\mu^\rm{F}_{\id,\vec{\bphi}}\\=
\int_{X^{k+1}}f_0\otimes (f_1\circ \tau^{t\bf{v}_{1,j}})\otimes
(f_2\circ \tau^{t\bf{v}_{2,j}})\otimes \cdots\otimes
(\sfE_\mu(f_k\,|\,\xi)\circ
\tau^{t\bf{v}_{k,j}})\,\d\mu^\rm{F}_{\id,\vec{\bphi}}
\end{multline*}
(because $\xi$ is respected by $\tau^{t\bf{v}_{k,j}}$ for each
$j$).  Therefore it suffices to prove the desired equality under the
assumption that $f_k$ is $\xi$-measurable.

However, with this in hand we may approximate $f$ as $g\cdot
\prod_{i=1}^{k-1}h_i$ in just the same was as for the proof of
Proposition~\ref{prop:AB-B}, and now rearranging according
to the different invariances of the functions $g$ and $h_i$ gives
\begin{eqnarray*}
&&\int_{X^{k+1}}f_0\otimes f_1\otimes f_2\otimes \cdots\otimes
(gh_1h_2\cdots h_{k-1})\,\d\mu^\rm{F}_{\id,\vec{\bphi}}\\
&&\quad\quad= \lim_{T\to\infty}\int_X f_0\cdot
A_{(0,T)}(f_1,f_2,\ldots,gh_1h_2\cdots h_{k-1})\,\d \mu\\
&&\quad\quad= \lim_{T\to\infty}\int_X f_0\cdot
A^-_{(0,T)}(f_1h_1,f_2h_2,\ldots,f_{k-1}h_{k-1},g)\,\d \mu\\
&&\quad\quad= \int_{X^{k+1}}f_0\otimes (f_1h_1)\otimes
(f_2h_2)\otimes\cdots\otimes
(f_{k-1}h_{k-1})\otimes g\,\d\mu^\rm{F}_{\id,\bphi_1,\bphi_2,\ldots,\bphi_{k-1},\bphi_k^-},
\end{eqnarray*}
where we let $A^-_{(0,T)}(f_1',f_2',\ldots,f_k')$ denote the multiple averages associated to the type-II precedent family
\[\{\bphi_1,\bphi_2,\ldots,\bphi_k^-\}.\]
(As in the previous proposition, we must tweak this step in case $\bphi_k^-$ is trivial, but then we simply let $A^-_{(0,T)}$ be associated to $\{\bphi_1,\bphi_2,\ldots,\bphi_{k-1}\}$, which is still a type-II precedent, and simply move $g$ outside these averages altogether.)

On the other hand, since
$h_i\circ \tau^{t\bf{v}_{k,j}} = h_i\circ \tau^{t\bf{v}_{i,j}}$, an analogous calculation gives
\begin{multline*}
\int_{X^{k+1}}f_0\otimes (f_1\circ \tau^{t\bf{v}_{1,j}})\otimes
\cdots\otimes
((gh_1h_2\cdots h_{k-1})\circ\tau^{t\bf{v}_{k,j}})\,\d\mu^\rm{F}_{\id,\vec{\bphi}}\\
= \int_{X^{k+1}}f_0\otimes ((f_1h_1)\circ
\tau^{t\bf{v}_{1,j}})\otimes\cdots\otimes
((f_{k-1}h_{k-1})\circ\tau^{t\bf{v}_{k-1,j}})\otimes (g\circ \tau^{t\bf{v}_{k,j}})\,\d\mu^\rm{F}_{\id,\bphi_1,\bphi_2,\ldots,\bphi_{k-1},\bphi_k^-}.
\end{multline*}

The proof is finished by treating two separate cases.

\quad\textbf{Case 1:}\quad$\underline{j \leq d-1}$\quad In this case the vectors $\bf{v}_{i,j}$ for $i\leq k$ are all still among the coefficient vectors of the precedent family $\{\bphi_1,\bphi_2,\ldots,\bphi_k^-\}$, and so applying the inductive assumption of conclusion $\rm{C}$ for that family gives the desired equality
\begin{multline*}
\int_{X^{k+1}}f_0\otimes (f_1h_1)\otimes
(f_2h_2)\otimes\cdots\otimes
(f_{k-1}h_{k-1})\otimes g\,\d\mu^\rm{F}_{\id,\bphi_1,\bphi_2,\ldots,\bphi_{k-1},\bphi_k^-}\\
= \int_{X^{k+1}}f_0\otimes ((f_1h_1)\circ
\tau^{t\bf{v}_{1,j}})\otimes\cdots\otimes
((f_{k-1}h_{k-1})\circ\tau^{t\bf{v}_{k-1,j}})\otimes (g\circ \tau^{t\bf{v}_{k,j}})\,\d\mu^\rm{F}_{\id,\bphi_1,\bphi_2,\ldots,\bphi_{k-1},\bphi_k^-}.
\end{multline*}

\quad\textbf{Case 2:}\quad$\underline{j = d}$\quad In this case we cannot make a simple inductive appeal to conclusion $\rm{C}$ for the precedent family, because the vector $\bf{v}_{k,d}$ is not among the coefficient vectors for that family: it has been removed from $\bphi_k^-$.  However, we also know that $g$ is $\tau^{\uhr\bbR\bf{v}_{k,d}}$-invariant, so in this case we write
\begin{multline*}
\int_{X^{k+1}}f_0\otimes ((f_1h_1)\circ
\tau^{t\bf{v}_{1,d}})\otimes\cdots\otimes
((f_{k-1}h_{k-1})\circ\tau^{t\bf{v}_{k-1,d}})\otimes (g\circ \tau^{t\bf{v}_{k,d}})\,\d\mu^\rm{F}_{\id,\bphi_1,\bphi_2,\ldots,\bphi_{k-1},\bphi_k^-}\\
= \int_{X^{k+1}}f_0\otimes ((f_1h_1)\circ
\tau^{t\bf{v}_{1,d}})\otimes\cdots\otimes
((f_{k-1}h_{k-1})\circ\tau^{t\bf{v}_{k-1,d}})\otimes g\,\d\mu^\rm{F}_{\id,\bphi_1,\bphi_2,\ldots,\bphi_{k-1},\bphi_k^-},
\end{multline*}
and now we can apply the inductive assumption of condition $\rm{C}$ to the precedent family, since all the remaining vectors $\bf{v}_{i,d}$, $i\leq k-1$, still appear as coefficients in that family.  This gives the desired equality again, and so completes the proof.  \qed

\textbf{Completion of the proof of Theorem~\ref{thm:B}}\quad  Suppose for the sake of contradiction that there were some good family of
f-polynomials $\vec{\bphi}$ for which
$\rm{A}(\vec{\bphi})\vee\rm{B}(\vec{\bphi})\vee\rm{C}(\vec{\bphi})$
fails.  Since $\prec$ satisfies the descending chain condition this
implies that there is a minimal such family, and by Lemma~\ref{lem:base-case} it must have at least two members.  However, for this family we know that
$\rm{A}(\prec\!\vec{\bphi})\vee\rm{B}(\prec\!\vec{\bphi})\vee\rm{C}(\prec\!\vec{\bphi})$
is true.  Therefore Proposition~\ref{prop:C-A} implies
that $\rm{A}(\vec{\bphi})$ is true; feeding this additional
information into Proposition~\ref{prop:AB-B} gives that
$\rm{B}(\vec{\bphi})$ is true; and now feeding this into
Propositino~\ref{prop:ABC-C} gives that
$\rm{C}(\vec{\bphi})$ is true.  This amounts to a contradiction,
and so completes the proof. \qed

\section{Further questions}

The ability to use fractional-power time-changes, $t\mapsto t^\a$, makes the study of polynomial multiple averages in higher-dimensional continuous flows much simpler than its discrete counterpart.  As remarked in the introduction, the older Conjecture~\ref{conj} for polynomial multiple averages in a $\bbZ^D$-action remains open, and it seems that the methods of the present paper offer little new insight into it.

However, within the setting of continuous actions there are other relatives of Theorem~\ref{thm:main} that may be of interest in their own right.

First, it is natural to ask whether the convergence of Theorem~\ref{thm:main} holds for arbitrary real intervals of increasing length: that is, whether the averages $A_I(f_1,\ldots,f_k)$ converge in $L^2(\mu)$ as $|I|\to\infty$, irrespective of how far the intervals $I$ lie from the origin in $\bbR$.  This question already creates difficulties for the use of fractional-power time-changes, since the correct replacement for Lemma~\ref{lem:change-vars} is not so clear outside the regime of tempered sequences of intervals.

A second question, posed by Vitaly Bergelson, asks whether the polynomial maps $\bf{p}_i$ appearing in that theorem could be replaced by members of some more general class.  Natural candidates here are members of the Hardy field of logarithmico-exponential functions, provided we choose such functions that do not grow too fast or too slowly.  More broadly still, it might be possible to single out a suitable class of maps $\bbR\to\bbR^D$ in terms of a Fej\'er-like condition on the asymptotic behaviour of their derivatives as $t\to\infty$.  In the setting of a single discrete transformation, some such extensions of the multiple recurrence and convergence theorems are known from the work of Bergelson and Knutson~\cite{BerKnu09}, who discuss the various growth conditions for these maps that are natural in this context.  We will not stop to formulate a precise conjecture here, as handling such more general maps would require too much preparation, but refer the reader to that paper for an introduction that can easily be transferred into the continuous-time world.

Another natural avenue for generalizing Theorem~\ref{thm:main} is to allow polynomial maps $\bbR^r\to \bbR^D$ for some higher-dimensional domain space $\bbR^r$.  In this case the formulation of our more general question is clear:

\begin{ques}
Is it true that for any action $\tau:\bbR^D\actson (X,\mu)$, any $f_1,f_2,\ldots,f_k \in L^\infty(\mu)$ and any polynomial maps $\bf{p}_i:\bbR^r\to\bbR^D$ for $i=1,2,\ldots,k$ that the associated multiple averages
\[\barint_{[0,R)^r}\prod_{i=1}^k(f_i\circ \tau^{\bf{p}_i(\bf{t})})\,\d t\]
converge in $L^2(\mu)$ as $R\to\infty$?
\end{ques}

I have not tried to extend the methods of the present paper to handle this question, but would not be surprised if that were possible.  Indeed, it should not be hard to extend the ideas of Section~\ref{sec:right-ques} to develop a notion of a `tempered' F\o lner sequence of domains $F_n\subset \bbR^r$, growing in size as $n\to\infty$ relatively fast compared to their distance from the origin, which is robust under a suitably rich collection of fractional-polynomial changes of variables
\[(t_1,t_2,\ldots,t_r)\mapsto (t_1^{\a_{1,1}}t_2^{\a_{1,2}}\cdots t_k^{\a_{1,k}},\ldots,t_1^{\a_{k,1}}t_2^{\a_{k,2}}\cdots t_k^{\a_{k,k}}).\]

However, even after introducing this flexibility, I do not at once see the right extension of the precedence ordering of Subsection~\ref{subs:prec} for fractional-polynomials in several variables, and so a little more work may be needed to find the right induction scheme if the above assertion to be proved this way.

Lastly, a more ambitious generalization in a similar vein would ask instead about an action of a connected nilpotent Lie group $\tau:N\actson (X,\mu)$ and polynomial maps $\bbR^r\to N$.  Starting with Leibman's work~\cite{Lei98} on multiple recurrence for actions of discrete nilpotent groups, there has been some progress on extending multiple recurrence and convergence theorems to nilpotent actions, and it might be hoped that the study of polynomial averages in an $\bbR^D$-action serves as a natural step towards such averages for a connected nilpotent group.  However, the greater complexity (albeit polynomial) of the multiplication in a nilpotent group makes it unclear whether the device of fractional-power time-changes yields similar simplifications as in the case of $\bbR^D$-actions, and I suspect that a more serious departure from the arguments above would be needed to settle this question.

\bibliographystyle{abbrv}
\bibliography{bibfile}

\end{document}